\documentclass[11pt]{amsart}

\usepackage{amssymb,amsfonts,amsmath}

\allowdisplaybreaks

\begin{document}
\theoremstyle{plain}
\newtheorem{thm}{Theorem}[section]
\newtheorem{theorem}[thm]{Theorem}
\newtheorem{lemma}[thm]{Lemma}
\newtheorem{corollary}[thm]{Corollary}
\newtheorem{proposition}[thm]{Proposition}
\newtheorem{conjecture}[thm]{Conjecture}
\theoremstyle{definition}
\newtheorem{construction}[thm]{Construction}
\newtheorem{notations}[thm]{Notations}
\newtheorem{question}[thm]{Question}
\newtheorem{problem}[thm]{Problem}
\newtheorem{remark}[thm]{Remark}
\newtheorem{remarks}[thm]{Remarks}
\newtheorem{definition}[thm]{Definition}
\newtheorem{claim}[thm]{Claim}
\newtheorem{assumption}[thm]{Assumption}
\newtheorem{assumptions}[thm]{Assumptions}
\newtheorem{properties}[thm]{Properties}
\newtheorem{example}[thm]{Example}
\newtheorem{comments}[thm]{Comments}
\newtheorem{blank}[thm]{}
\newtheorem{observation}[thm]{Observation}
\newtheorem{defn-thm}[thm]{Definition-Theorem}

\newcommand{\sM}{{\mathcal M}}


\title[Mirzakharni's recursion formula]{Mirzakharni's recursion formula is
equivalent to the Witten-Kontsevich theorem}
        \author{Kefeng Liu}
        \address{Center of Mathematical Sciences, Zhejiang University, Hangzhou, Zhejiang 310027, China;
                Department of Mathematics,University of California at Los Angeles,
                Los Angeles, CA 90095-1555, USA}
        \email{liu@math.ucla.edu, liu@cms.zju.edu.cn}
        \author{Hao Xu}
        \address{Center of Mathematical Sciences, Zhejiang University, Hangzhou, Zhejiang 310027, China}
        \email{haoxu@cms.zju.edu.cn}

        \dedicatory{Dedicated to Jean-Michel Bismut on the occasion
of his 60th birthday}

        \begin{abstract}
        In this paper, we give a proof of Mirzakhani's recursion formula of Weil-Petersson volumes of moduli spaces
        of curves using the Witten-Kontsevich theorem.
        We also describe properties of intersections numbers involving higher degree $\kappa$ classes.
        \end{abstract}


    \maketitle

\section{Introduction}

Following the notation of Mulase and Safnuk \cite{MS}, let $\mathcal
M_{g,n}(\bold L)$ denote the moduli space of bordered Riemann
surfaces with $n$ geodesic boundary components of specified lengths
$\bold L=(L_1,\dots,L_n)$ and let ${\rm Vol}_{g,n}(\bold L)$ denote
its Weil-Petersson volume ${\rm Vol}(\mathcal M_{g,n}(\bold L))$.
Using her remarkable generalization of the McShane identity,
Mirzakhani \cite{Mir} proved a beautiful recursion formula for these
Weil-Petersson volumes
\begin{multline*}
    {\rm Vol}_{g,n}(\bold L)
    = \frac{1}{2 L_1} \sum_{\substack{g_1+g_2 = g \\ \underline{n}=I \coprod J } }
        \int_{0}^{L_1} \int_{0}^{\infty} \int_{0}^{\infty} xy H(t, x+y) \\
        \times{\rm Vol}_{g_1,n_1 }(x, \bold L_{I} )
        {{\rm Vol}}_{g_2,n_2 }(y, \bold L_{J}) dxdydt
    \\    +\frac{1}{2L_1} \int_{0}^{L_1} \int_{0}^{\infty} \int_{0}^{\infty} xy H(t, x+y) {{\rm Vol}}_{g-1,n+1}(x,y,L_2,\dots,L_n)dxdydt
    \\
    +\frac{1}{2L_1} \sum_{j=2}^{n} \int_{0}^{L_1}\int_{0}^{\infty}
        x \bigl( H(x, L_1+L_j) + H(x, L_1 - L_j) \bigr) \\
        \times {{\rm Vol}}_{g,n-1}(x, L_2,\dots,\hat{L_j},\dots,L_n) dxdt,
\end{multline*}
where the kernel function
$$H(x,y) = \frac{1 }{1 + e^{(x+y)/2 } }  +  \frac{1 }{1 + e^{(x-y)/2 }}.$$

Using symplectic reduction, Mirzakhani \cite{Mir2} showed the
following relation
\begin{align*}
\frac{{\rm Vol}_{g,n}(2\pi\bold L)}{(2\pi^2)^{3g+n-3} }
  &= \frac{1}{(3g+n-3)!} \int_{\mathcal M_{g,n}} ( \kappa_1 + \sum_{i=1}^n L_i^2 \psi_i
    )^{3g+n-3}\\
    &= \sum_{\substack{d_0 + \cdots + d_n \\= 3g+n-3}} \prod_{i=0}^{n}\frac{1}{d_i! }
        \langle\kappa_1^{d_0} \prod\tau_{d_i}\rangle_{g,n} \prod_{i=1}^{\infty}
      L_{i}^{2d_i}.
\end{align*}
Combining with her recursion formula of Weil-Petersson volumes,
Mirzakhani \cite{Mir2} found a new proof of the celebrated
Witten-Kontsevich theorem.

By taking derivatives with respect to $\bold L=(L_1,\dots,L_n)$ in
Mirzakhani's recursion, Mulase and Safnuk \cite{MS} obtained the
following enlightening recursion formula of intersection numbers
which is equivalent to Mirzakhani's recursion.
\begin{multline*}
 (2d_1+1)!!\langle\prod_{j=1}^n\tau_{d_j}\kappa_1^{a}\rangle_g \\
  =  \sum_{j=2}^{n} \sum_{b=0}^{a}\frac{a!}{(a-b)!} \frac{(2(b+d_1+d_j)-1)!!}{(2d_j-1)!!} \beta_b \langle\kappa_1^{a-b}\tau_{b+d_1+d_j-1}\prod_{i \neq 1,j}\tau_{d_i}\rangle_g \\
 + \frac{1}{2} \sum_{b=0}^a \sum_{r + s = b+d_1 - 2}  \frac{a!}{(a-b)!}
  (2r+1)!! (2s+1)!! \beta_b \langle\kappa_1^{a-b}\tau_r\tau_s\prod_{i\neq 1}\tau_{d_i}\rangle_{g-1}  \\
  + \frac{1}{2}\sum_{b=0}^a\sum_{\substack{c+c' = a-b \\ I \coprod J
= \{2, \ldots, n \}}}
  \sum_{r + s = b+d_1-2}  \frac{a!}{c! c'!}
   (2r+1)!!(2s+1)!! \beta_b \\
  \times \langle\kappa_1^{c}\tau_r\prod_{i\in I} \tau_{d_i}\rangle_{g'}
 \langle\kappa_1^{c'}\tau_s\prod_{i\in J} \tau_{d_i}\rangle_{g-g'},
\end{multline*}
where
\begin{equation*}
\beta_b = (2^{2b+1}-4)\frac{\zeta(2b)}{(2\pi^2)^b} =
(-1)^{b-1}2^b(2^{2b}-2) \frac{B_{2b}}{(2b)!}.
\end{equation*}

Safnuk \cite{Saf} gave a proof of the above differential form of
Mirzakhani's recurson formula using localization techniques, but he
also used the Mirzakhani-McShane formula. The relationship between
Mirzakhani's recurson and matrix integrals has been studied by
Eynard-Orantin \cite{EO} and Eynard \cite{Ey}.

Indeed, when $a=0$, Mulase-Safnuk differential form of Mirzakhani's
recursion is just the  Witten-Kontsevich theorem \cite{Wi, Ko} in
the form of DVV recursion relation \cite{DVV}. There are several
other new proofs of Witten-Kontsevich theorem \cite{CLL, KaL, KL,
OP} besides Mirzakhani's proof \cite{Mir2}.

More discussions about Weil-Petersson volumes from the point of view
of intersection numbers can be found in the papers \cite{DN, KK, MZ,
Zo}.

In Section 2, we show that Mirzakhani's recursion formula is
essentially equivalent to the Witten-Kontsevich theorem via a
formula from \cite{KMZ} expressing $\kappa$ classes in terms of
$\psi$ classes. In Section 3, we present certain results of
intersection numbers involving higher
degree $\kappa$ classes.\\

\noindent{\bf Acknowledgements.} We would like to thank Chiu-Chu
Melissa Liu for helpful discussions. We also thank the referees for
helpful suggestions.

\section{Proof of Mirzakhani's recursion formula}
We first give three lemmas. The following lemma can be found in
\cite{MS}.
\begin{lemma}
The constants $\beta_b$ in Mirzakhani's recursion satisfy the
following
$$\sum_{k=0}^\infty \beta_k x^k=\frac{\sqrt{2x}}{\sin \sqrt{2x}}.$$
And its inverse
$$(\sum_{k=0}^\infty \beta_k x^k)^{-1}=\frac{\sin \sqrt{2x}}{\sqrt{2x}}=\sum_{k=0}^\infty\frac{(-1)^k 2^k}{(2k+1)!}x^k$$
\end{lemma}
\begin{proof}
Since
$$    \sum_{n=0}^{\infty} \frac{B_{2n}}{(2n)! } x^{2n}
    = \frac{x}{2}\frac{e^{x/2} + e^{-x/2} }{e^{x/2} - e^{-x/2}
    }=\frac{x}{2i}\cot\frac{x}{2i},$$
we have
$$\sum_{k=0}^\infty\beta_k x^k=\sqrt{2x}(\cot\sqrt{\frac{x}{2}}-\cot\sqrt{2x})=\frac{\sqrt{2x}}{\sin\sqrt{2x}}.$$
\end{proof}

The following elementary result is crucial to our proof.
\begin{lemma}
Let $F(m,n)$ and $G(m,n)$ be two functions defined on $\mathbb
N\times\mathbb N$, where $\mathbb N=\{0,1,2,\dots\}$ is the set of
nonnegative integers. Let $\alpha_k$ and $\beta_k$ be real numbers
that satisfy
$$\sum_{k=0}^\infty\alpha_k x^k=(\sum_{k=0}^\infty\beta_k x^k)^{-1}.$$
Then the following two identities are equivalent.

\begin{align*}
G(m,n)=\sum_{k=0}^m \alpha_k F(m-k,n+k),\quad \forall\ (m,n)\in\mathbb N\times\mathbb N\\
F(m,n)=\sum_{k=0}^m \beta_k G(m-k,n+k),\quad \forall\ (m,n)\in\mathbb N\times\mathbb N\\
\end{align*}
\end{lemma}
\begin{proof}
Assume the first identity holds, then we have
\begin{align*}
\sum_{i=0}^m \beta_i G(m-i,n+i)&=\sum_{i=0}^m \beta_i \sum_{j=0}^{m-i}\alpha_j F(m-i-j,n+i+j)\\
&=\sum_{k=0}^m\sum_{i+j=k}(\beta_i\alpha_j) F(m-k,n+k)\\
&=\sum_{k=0}^m\delta_{k0}F(m-k,n+k)\\
&=F(m,n).
\end{align*}
So we proved the second identity. The proof of the other direction
is the same.
\end{proof}

The fact that intersection numbers involving both $\kappa$ classes
and $\psi$ classes can be reduced to intersection numbers involving
only $\psi$ classes was already known to Witten [9], and has been
developed by Arbarello-Cornalba \cite{Ar-Co}, Faber \cite{Fa} and
Kaufmann-Manin-Zagier \cite{KMZ} into a nice combinatorial
formalism.

\begin{lemma} \cite{KMZ} For $m>0$,
$$\langle\prod_{j=1}^n\tau_{d_j}\kappa_1^m\rangle_g=\sum_{k=1}^m\frac{(-1)^{m-k}}{k!}\sum_{\substack {m_1+\cdots+m_k=m\\m_i>0}}\binom{
m}{m_1,\dots,m_k}\langle\prod_{j=1}^n\tau_{d_j}\prod_{j=1}^k\tau_{m_j+1}\rangle_g.$$
\end{lemma}
\begin{proof} (sketch) Let
$ \pi_{n+p,n}:\overline{\sM}_{g,n+p}\longrightarrow
\overline{\sM}_{g,n} $ be the morphism which forgets the last $p$
marked points and denote
$\pi_{n+p,n*}(\psi_{n+1}^{a_1+1}\dots\psi_{n+p}^{a_p+1})$ by
$R(a_1,\dots,a_p)$, then we have the formula from \cite{Ar-Co}
$$R(a_1,\dots,a_p)=\sum_{\sigma\in\mathbb S_p}\prod_{\substack{{\rm each\ cycle}\ c\\ {\rm of}\ \sigma}}\kappa_{\sum_{j\in c}a_j},$$
where we write any permutation $\sigma$ in the symmetric group
$\mathbb S_p$ as a product of disjoint cycles.

A formal combinatorial argument \cite{KMZ} leads to the following
inversion equation
$$\kappa_{a_1}\cdots\kappa_{a_p}=\sum_{k=1}^p\frac{(-1)^{p-k}}{k!}\sum_{\substack{\{1,\dots,p\}=S_1\coprod\dots\coprod S_k\\S_k\neq\emptyset}}
R(\sum_{j\in S_1}a_j,\dots,\sum_{j\in S_k}a_j),$$ from which the
result follows easily.
\end{proof}

\begin{proposition}
\begin{multline*}
 \sum_{b=0}^a(-1)^b\binom{a}{b}\frac{(2(d_1+b)+1)!!}{(2b+1)!!}\langle\tau_{d_1+b}\prod_{i=2}^n\tau_{d_i}\kappa_1^{a-b}\rangle_g\\
  =  \sum_{j=2}^{n} \frac{(2d_1+2d_j-1)!!}{(2d_j-1)!!} \langle\kappa_1^{a}\tau_{d_1+d_j-1}\prod_{i \neq 1,j}\tau_{d_i}\rangle_g \\
 + \frac{1}{2} \sum_{r + s = d_1 - 2}
  (2r+1)!! (2s+1)!! \langle\kappa_1^{a}\tau_r\tau_s\prod_{i\neq 1}\tau_{d_i}\rangle_{g-1}  \\
  + \frac{1}{2}\sum_{\substack{ c+c'=a\\I \coprod J= \{2, \ldots, n
  \}}}\binom{a}{c}
  \sum_{r + s = d_1-2} (2r+1)!!(2s+1)!!
\langle\kappa_1^{c}\tau_r\prod_{i\in I} \tau_{d_i}\rangle_{g'}
 \langle\kappa_1^{c'}\tau_s\prod_{i\in J} \tau_{d_i}\rangle_{g-g'}.
\end{multline*}
\end{proposition}
\begin{proof}
Let LHS and RHS denote the left and right hand side of the equation
respectively. By Lemma 2.3 and the Witten-Kontsevich theorem, we
have
\begin{multline*}
(2d_1+1)!!\langle\prod_{j=1}^n\tau_{d_j}\kappa_1^{a}\rangle_g\\
=(2d_1+1)!!\sum_{k=0}^a\frac{(-1)^{a-k}}{k!}\sum_{\substack
{m_1+\cdots+m_k=a\\m_i>0}}\binom{
a}{m_1,\dots,m_k}\langle\prod_{j=1}^n\tau_{d_j}\prod_{j=1}^k\tau_{m_j+1}\rangle_g\\
=\sum_{k=0}^a \frac{(-1)^{a-k}}{k!}\sum_{\substack
{m_1+\cdots+m_k=a\\m_i>0}}\binom{
a}{m_1,\dots,m_k}\\\times\left(\sum_{j=2}^n\frac{(2(d_1+d_j)-1)!!}{(2d_j-1)!!}
\langle\tau_{d_1+d_j-1}\prod_{i\neq
1,j}\tau_{d_i}\prod_{i=1}^k\tau_{m_i+1}\rangle_g\right.\\
+\sum_{j=1}^k\frac{(2(d_1+m_j)+1)!!}{(2m_j+1)!!}\langle\tau_{d_1+m_j}\prod_{i=2}^n\tau_{d_i}\prod_{i\neq
j}\tau_{m_i+1}\rangle_g\\
+\frac{1}{2}\sum_{r+s=d_1-2}(2r+1)!!(2s+1)!!\langle\tau_r\tau_s\prod_{i=2}^n\tau_{d_i}\prod_{i=1}^k\tau_{m_i+1}\rangle_{g-1}\\
\left.+\frac{1}{2}\sum_{\substack{I\coprod
J=\{2,\dots,n\}\\I'\coprod J'=\{1,\dots,k\}
}}\sum_{r+s=d_1-2}(2r+1)!!(2s+1)!!\right.\\
\left.\times\langle\tau_r\prod_{i\in I}\tau_{d_i}\prod_{i\in
I'}\tau_{m_i+1}\rangle_{g'}\langle\tau_s\prod_{i\in
J}\tau_{d_i}\prod_{i\in J'}\tau_{m_i+1}\rangle_{g-g'}\right)\\
=\sum_{j=2}^{n} \frac{(2d_1+2d_j-1)!!}{(2d_j-1)!!} \langle\kappa_1^{a}\tau_{d_1+d_j-1}\prod_{i \neq 1,j}\tau_{d_i}\rangle_g \\
 + \frac{1}{2} \sum_{r + s = d_1 - 2}
  (2r+1)!! (2s+1)!! \langle\kappa_1^{a}\tau_r\tau_s\prod_{i\neq 1}\tau_{d_i}\rangle_{g-1}  \\
  + \frac{1}{2}\sum_{\substack{ c+c'=a\\I \coprod J= \{2, \ldots, n
  \}}}\binom{a}{c}
  \sum_{r + s = d_1-2} (2r+1)!!(2s+1)!!
\langle\kappa_1^{c}\tau_r\prod_{i\in I} \tau_{d_i}\rangle_{g'}
 \langle\kappa_1^{c'}\tau_s\prod_{i\in J} \tau_{d_i}\rangle_{g-g'}\\
+\sum_{k=0}^a \frac{(-1)^{a-k}}{k!}\sum_{\substack
{m_1+\cdots+m_k=a\\m_i>0}}\binom{
a}{m_1,\dots,m_k}\\
\times\sum_{j=1}^k\frac{(2(d_1+m_j)+1)!!}{(2m_j+1)!!}\langle\tau_{d_1+m_j}\prod_{i=2}^n\tau_{d_i}\prod_{i\neq
j}\tau_{m_i+1}\rangle_g\\
 =RHS+\sum_{k\geq0} \frac{(-1)^{a-k-1}}{(k+1)!}\sum_{b=1}^a\sum_{\substack
{m_1+\cdots+m_k=a-b\\m_i>0}}\binom{a}{b}\binom{
a-b}{m_1,\dots,m_k}\\
\times(k+1)\frac{(2(d_1+b)+1)!!}{(2b+1)!!}\langle\tau_{d_1+b}\prod_{i=2}^n\tau_{d_i}\prod_{i=1
}^k\tau_{m_i+1}\rangle_g\\
=RHS-\sum_{b=1}^a(-1)^b\binom{a}{b}\frac{(2(d_1+b)+1)!!}{(2b+1)!!}\langle\tau_{d_1+b}\prod_{i=2}^n\tau_{d_i}\kappa_1^{a-b}\rangle_g\\
=RHS-LHS+(2d_1+1)!!\langle\prod_{j=1}^n\tau_{d_j}\kappa_1^{a}\rangle_g.
\end{multline*}

So we have proved $RHS=LHS$.
\end{proof}

Proposition 2.4 is also implicitly contained in the arguments of
Mulase and Safnuk \cite{MS}.

\begin{theorem}
\begin{multline*}
 \frac{(2d_1+1)!!}{a!}\langle\prod_{j=1}^n\tau_{d_j}\kappa_1^{a}\rangle_g \\
  =  \sum_{b=0}^{a}\sum_{j=2}^{n}  \frac{(2(b+d_1+d_j)-1)!!}{(a-b)!(2d_j-1)!!} \beta_b \langle\kappa_1^{a-b}\tau_{b+d_1+d_j-1}\prod_{i \neq 1,j}\tau_{d_i}\rangle_g \\
 + \frac{1}{2} \sum_{b=0}^a \sum_{r + s = b+d_1 - 2}  \frac{(2r+1)!! (2s+1)!!}{(a-b)!}
   \beta_b \langle\kappa_1^{a-b}\tau_r\tau_s\prod_{i\neq 1}\tau_{d_i}\rangle_{g-1}  \\
  + \frac{1}{2}\sum_{b=0}^a\sum_{\substack{c+c' = a-b \\ I \coprod J
= \{2, \ldots, n \}}}
  \sum_{r + s = b+d_1-2}  \frac{(2r+1)!!(2s+1)!!}{c! c'!}
   \beta_b \\
  \times \langle\kappa_1^{c}\tau_r\prod_{i\in I} \tau_{d_i}\rangle_{g'}
 \langle\kappa_1^{c'}\tau_s\prod_{i\in J} \tau_{d_i}\rangle_{g-g'},
\end{multline*}
where the constants $\beta_k$ are given by
$$(\sum_{k=0}^\infty \beta_k x^k)^{-1}=\frac{\sin \sqrt{2x}}{\sqrt{2x}}=\sum_{k=0}^\infty\frac{(-1)^k}{k!(2k+1)!!}x^k.$$
\end{theorem}
\begin{proof}
Denote the LHS by $F(a,d_1)$. Let
\begin{multline*}
 G(a,d_1)=
\sum_{j=2}^{n}  \frac{(2(d_1+d_j)-1)!!}{a!(2d_j-1)!!} \langle\kappa_1^{a}\tau_{d_1+d_j-1}\prod_{i \neq 1,j}\tau_{d_i}\rangle_g \\
 + \frac{1}{2} \sum_{r + s = d_1 - 2}  \frac{(2r+1)!! (2s+1)!!}{a!}
 \langle\kappa_1^{a}\tau_r\tau_s\prod_{i\neq 1}\tau_{d_i}\rangle_{g-1}  \\
  + \frac{1}{2}\sum_{\substack{c+c' = a \\ I \coprod J
= \{2, \ldots, n \}}}
  \sum_{r + s = d_1-2}  \frac{(2r+1)!!(2s+1)!!}{c! c'!}
  \times \langle\kappa_1^{c}\tau_r\prod_{i\in I} \tau_{d_i}\rangle_{g'}
 \langle\kappa_1^{c'}\tau_s\prod_{i\in J} \tau_{d_i}\rangle_{g-g'},
\end{multline*}

Note that Proposition 2.4 is just
$$\sum_{b=0}^a\frac{(-1)^b}{b!(2b+1)!!}F(a-b,d_1+b)=G(a,d_1).$$
By Lemmas 2.1 and 2.2, we have
$$F(a,d_1)=\sum_{b=0}^a\beta_b G(a-b,d_1+b)=RHS.$$
So we conclude the proof.
\end{proof}

\section{Higher Weil-Petersson volumes}

Mirzakhani's formula provides a recursive way of computing the
following Weil-Petersson volumes of moduli spaces of curves
$$WP(g):=\int_{\overline{\sM}_{g,n}}\kappa_1^{3g-3+n}.$$
Mirzakhani's formula resorts to intersection numbers of mixed $\psi$
and $\kappa$ classes.

A natural question is whether there exist an explicit formula
expressing $WP(g)$ in terms of those $WP(g')$ with $g'<g$. Recall
the following beautiful formula due to Itzykson-Zuber \cite{IZ}.

\begin{proposition}{\bf (Itzykson-Zuber)} Let $g\geq0$. Then
$$\phi_{g+1}=\frac{25g^2-1}{24}\phi_g+\frac{1}{2}\sum_{m=1}^g\phi_{g+1-m}\phi_m,$$
where $\phi_0=-1, \phi_1=\frac{1}{24}$ and
$$\phi_g=\frac{(5g-5)(5g-3)}{2^g(3g-3)!}\langle\tau_2^{3g-3}\rangle_g,
\quad g\geq2.$$
\end{proposition}
By projection formula, we have
$$\langle\tau_2^{3g-3}\rangle_g=\langle\kappa_1^{3g-3}\rangle_g+\cdots,$$
where $\cdots$ denote terms involving higher degree kappa classes.
Also note that $\langle\kappa_1^{3g-3}\rangle_g$ is conjecturally
\cite{LX2} the largest term in the right hand side.

To our disappointment, so far, all recursion formulae for $WP(g)$
stemming from the Witten-Kontsevich theorem involve either $\psi$
class or higher degree $\kappa$ classes inevitably.

Mirzakhani, Mulase and Safnuk's arguments use Wolpert's formula
\cite{Wo}
$$\kappa_1=\frac{1}{2\pi^2}\omega_{WP},$$
where $\omega_{WP}$ is the Weil-Petersson K\"ahler form. We have no
similar formulae for higher degree $\kappa$ classes. So a priori
$\kappa_1$ may be rather special in the intersection theory.
However, as we will see, this is not the case.

First we fix notations as in \cite{KMZ}. Consider the semigroup
$N^\infty$ of sequences ${\bold m}=(m(1),m(2),\dots)$ where $m(i)$
are nonnegative integers and $m(i)=0$ for sufficiently large $i$.

Let $\bold m, \bold t, \bold{a_1,\dots,a_n} \in N^\infty$, $\bold
m=\sum_{i=1}^n \bold{a_i}$, $\bold m\geq\bold t$ and $\bold
s:=(s_1,s_2,\dots)$ be a family of independent formal variables.
$$|\bold m|:=\sum_{i\geq 1}i m(i),\quad ||\bold m||:=\sum_{i\geq1}m(i),\quad \bold s^{\bold m}:=\prod_{i\geq 1}s_i^{m(i)},\quad \bold m!:=\prod_{i\geq1}m(i)!,$$
$$\binom{\bold m}{\bold{t}}:=\prod_{i\geq1}\binom{ m(i)}{t(i)},\quad \binom{\bold m}{\bold{a_1,\dots,a_n}}:=\prod_{i\geq1}\binom{ m(i)}{a_1(i),\dots,a_n(i)}.$$

Let $\bold b\in N^\infty$, we denote a formal monomial of $\kappa$
classes by
$$\kappa(\bold b):=\prod_{i\geq1}\kappa_i^{b(i)}.$$

We are interested in the following intersection numbers
$$\langle\kappa(\bold b)\tau_{d_1}\cdots\tau_{d_n}\rangle_g:=\int_{\overline{\sM}_{g,n}}\kappa(\bold b)\psi_1^{d_1}\cdots\psi_n^{d_n}.$$

When $d_1=\cdots =d_n=0$, these intersection numbers are called
higher Weil-Petersson volumes of moduli spaces of curves. The
details of the following discussions are contained in \cite{LX}.

The following lemma is a direct generalization of Lemma 2.2.
\begin{lemma}
Let $F(\bold L,n)$ and $G(\bold L,n)$ be two functions defined on
$N^\infty\times\mathbb N$, where $\mathbb N=\{0,1,2,\dots\}$ is the
set of nonnegative integers. Let $\alpha_{\bold L}$ and
$\beta_{\bold L}$ be real numbers depending only on $\bold L\in
N^\infty$ that satisfy $\alpha_{\bold 0}\beta_{\bold 0}=1$ and
$$\sum_{\bold L+\bold{L'}=\bold b}\alpha_{\bold L}\beta_{\bold{L'}}=0,\qquad\bold b\neq0.$$
Then the following two identities are equivalent.

\begin{align*}
G(\bold b,n)=\sum_{\bold L+\bold{L'}=\bold b} \alpha_{\bold L} F(\bold{L'},n+|\bold L|),\quad\forall\ (\bold b,n)\in N^\infty\times\mathbb N\\
F(\bold b,n)=\sum_{\bold L+\bold{L'}=\bold b} \beta_{\bold L} G(\bold{L'},n+|\bold L|),\quad\forall\ (\bold b,n)\in N^\infty\times\mathbb N\\
\end{align*}
\end{lemma}

We may generalize Mirzakhani's recursion formula to include higher
degree $\kappa$ classes.
\begin{theorem} There exist (uniquely determined) rational numbers $\alpha_{\bold
L}$ depending only on $\bold L\in N^\infty$, such that for any
$\bold b\in N^\infty$ and $d_j\geq 0$,  the following recursion
relation of mixed $\psi$ and $\kappa$ intersection numbers holds.
\begin{multline*}
(2d_1+1)!!\langle\kappa(\bold b)\prod_{j=1}^n\tau_{d_j}\rangle_g\\
=\sum_{j=2}^n\sum_{\bold L+\bold{L'}=\bold b}\alpha_{\bold
L}\binom{\bold b}{\bold L}\frac{(2(|\bold
L|+d_1+d_j)-1)!!}{(2d_j-1)!!}\langle\kappa(\bold{L'})\tau_{|\bold
L|+d_1+d_j-1}\prod_{i\neq
1,j}\tau_{d_i}\rangle_g\\
+\frac{1}{2}\sum_{\bold L+\bold{L'}=\bold b}\sum_{r+s=|\bold
L|+d_1-2}\alpha_{\bold
L}\binom{\bold b}{\bold L}(2r+1)!!(2s+1)!!\langle\kappa(\bold{L'})\tau_r\tau_s\prod_{i\neq1}\tau_{d_i}\rangle_{g-1}\\
+\frac{1}{2}\sum_{\substack{\bold L+\bold{e}+\bold{f}=\bold
b\\I\coprod J=\{2,\dots,n\}}}\sum_{r+s=|\bold L|+d_1-2}\alpha_{\bold
L}\binom{\bold b}{\bold
L,\bold{e},\bold{f}}(2r+1)!!(2s+1)!!\\
\times \langle\kappa(\bold{e})\tau_r\prod_{i\in
I}\tau_{d_i}\rangle_{g'}\langle\kappa(\bold{f})\tau_s\prod_{i\in
J}\tau_{d_i}\rangle_{g-g'}.
\end{multline*}
These tautological constants $\alpha_{\bold L}$ can be determined
recursively from the following formula
$$\sum_{\bold L+\bold{L'}=\bold b}\frac{(-1)^{||\bold L||}\alpha_{\bold L}}{\bold L!\bold{L'}!(2|\bold{L'}|+1)!!}=0,\qquad \bold b\neq0,$$
namely
$$\alpha_{\bold b}=\bold b!\sum_{\substack{\bold L+\bold{L'}=\bold b\\ \bold{L'}\neq\bold 0}}\frac{(-1)^{||\bold L'||-1}\alpha_{\bold L}}{\bold L!\bold{L'}!(2|\bold{L'}|+1)!!},\qquad\bold b\neq0,$$
with the initial value $\alpha_{\bold 0}=1$.
\end{theorem}

\begin{theorem}
\begin{multline*}
\sum_{\bold L+\bold{L'}=\bold b}(-1)^{||\bold L||}\binom{\bold b}{\bold L}\frac{(2d_1+2|\bold L|+1)!!}{(2|\bold L|+1)!!}
\langle\kappa(\bold L')\tau_{d_1+|\bold L|}\prod_{j=2}^n\tau_{d_j}\rangle_g\\
=\sum_{j=2}^n\frac{(2(d_1+d_j)-1)!!}{(2d_j-1)!!}\langle\kappa(\bold{b})\tau_{d_1+d_j-1}\prod_{i\neq
1,j}\tau_{d_i}\rangle_g\\
+\frac{1}{2}\sum_{r+s=|d_1|-2}(2r+1)!!(2s+1)!!\langle\kappa(\bold{b})\tau_r\tau_s\prod_{i\neq1}\tau_{d_i}\rangle_{g-1}\\
+\frac{1}{2}\sum_{\substack{\bold{e}+\bold{f}=\bold b\\I\coprod
J=\{2,\dots,n\}}}\sum_{r+s=d_1-2}\binom{\bold b}{\bold{e}}(2r+1)!!(2s+1)!!\\
\times \langle\kappa(\bold{e})\tau_r\prod_{i\in
I}\tau_{d_i}\rangle_{g'}\langle\kappa(\bold{f})\tau_s\prod_{i\in
J}\tau_{d_i}\rangle_{g-g'}.
\end{multline*}
\end{theorem}

Theorem 3.3 and Theorem 3.4 implies each other through Lemma 3.2.

Both Theorems 3.3 and 3.4 are effective recursion formulae for
computing higher Weil-Petersson volumes with the three initial
values
$$\langle\tau_0\kappa_1\rangle_1=\frac{1}{24},\qquad
\langle\tau_0^3\rangle_0=1,\qquad
\langle\tau_1\rangle_1=\frac{1}{24}.$$ From the following
Proposition 3.4, we have
$$\langle\kappa(\bold
b)\rangle_g=\frac{1}{2g-2}\sum_{\bold L+\bold L'=\bold
b}(-1)^{||\bold L||}\binom{\bold b}{\bold L}\langle\tau_{|\bold
L|+1}\kappa(\bold L')\rangle_g.$$

We have computed a table of $\alpha_{\bold L}$ for all $|\bold
L|\leq 15$ and have written a Maple program \cite{Maple}
implementing Theorems 3.3 and 3.4.

In fact, we find that $\psi$ and $\kappa$ classes are compatible in
the sense that recursions of pure $\psi$ classes can be neatly
generalized to recursions including both $\psi$ and $\kappa$ classes
by the same proof as Proposition 2.4. In view of Theorem 3.8 below,
this can be rephrased as differential equations governing generating
functions of $\psi$ classes also govern generating functions of
mixed $\psi$ and $\kappa$ classes.

We present some examples below.

\begin{proposition}
Let $\bold b\in N^\infty$ and $d_j\geq0$. Then
$$\sum_{\bold L+\bold L'=\bold b}(-1)^{||\bold
L||}\binom{\bold b}{\bold L}\langle\tau_{|\bold
L|+1}\prod_{j=1}^n\tau_{d_j}\kappa(\bold
L')\rangle_g=(2g-2+n)\langle\prod_{j=1}^n\tau_{d_j}\kappa(\bold
b)\rangle_g.$$
\end{proposition}

The above proposition is a generalization of the dilaton equation.
In the special case $\bold b=(m,0,0,\dots)$, it has been proved by
Norman Do and Norbury \cite{DN}.

\begin{proposition}
Let $\bold b\in N^\infty$. Then
\begin{multline*}
\langle\tau_0\tau_1\prod_{j=1}^n\tau_{d_j}\kappa(\bold{
b})\rangle_g=\frac{1}{12}\langle\tau_0^4\prod_{
j=1}^n\tau_{d_j}\kappa(\bold b)\rangle_g\\
+\frac{1}{2}\sum_{\substack{\bold L+\bold{L'}=\bold
b\\\underline{n}=I\coprod J}}\binom{\bold b}{\bold
L}\langle\tau_0^2\prod_{i\in I}\tau_{d_i}\kappa(\bold
L)\rangle_{g'}\langle\tau_0^2\prod_{i\in
J}\tau_{d_i}\kappa(\bold{L'})\rangle_{g-g'}.
\end{multline*}
\end{proposition}

The above proposition, together with the projection formula, can be
used to derive an effective recursion formula for higher
Weil-Petersson volumes \cite{LX} (without $\psi$ classes).

Let $\bold s:=(s_1,s_2,\dots)$ and $\bold t:=(t_0,t_1,t_2,\dots)$,
we introduce the following generating function
$$G(\bold s,\bold t):=\sum_{g}\sum_{\bold m,\bold n}\langle\kappa_1^{m_1}\kappa_2^{m_2}\cdots\tau_0^{n_0}\tau_1^{n_1}\cdots\rangle_g\frac{\bold s^{\bold m}}{\bold m!}\prod_{i=0}^\infty\frac{t_i^{n_i}}{n_i!},$$
where $\bold s^{\bold m}=\prod_{i\geq1}s_i^{m_i}$.

Following Mulase and Safnuk \cite{MS}, we introduce the following
family of differential operators for $k\geq -1$,
 \begin{multline*}
    V_k = -\frac{1}{2} \sum_{\bold L} (2(|\bold L|+k)+3)!! \frac{(-1)^{||\bold L||}}{\bold L!(2|\bold
 L|+1)!!} \bold
    s^{\bold L}
        \frac{\partial }{\partial t_{|\bold L|+k+1} }\\
        + \frac{1}{2} \sum_{j=0}^{\infty} \frac{(2(j+k)+1)!! }{(2j-1)!! } t_j
        \frac{\partial }{\partial t_{j+k} }
    + \frac{1}{4} \sum_{d_1 + d_2 = k-1}
        (2d_1 + 1)!! (2d_2 + 1)!! \frac{\partial^2 }{\partial t_{d_1} \partial t_{d_2}}\\
        + \frac{\delta_{k,-1}t_0^2}{4} + \frac{\delta_{k,0} }{48}.
 \end{multline*}

\begin{theorem} \cite{LX, MS}
The recursion of Theorem 3.4 implies
$$V_k\exp(G)=0.$$ Moreover, we can check directly that the operators
$V_k$, $k\geq-1$ satisfy the Virasoro relations
$$[V_n,V_m]=(n-m)V_{n+m}.$$
\end{theorem}

The Witten-Kontsevich theorem
 states that the generating function for $\psi$ class
intersections
$$F(t_0, t_1, \ldots)= \sum_{g} \sum_{\bold n} \langle\prod_{i=0}^\infty \tau_{i}^{n_i}\rangle_{g} \prod_{i=0}^\infty \frac{t_i^{n_i} }{n_i!
        }$$
is a $\tau$-function for the KdV hierarchy.

\begin{theorem} \cite{LX, MS} We have
$$G(\bold s,t_0,t_1,\dots)=F(t_0,t_1,t_2+p_2,t_3+p_3,\dots),$$
where $p_k$ are polynomials in $\bold s$ given by
$$p_k=\sum_{|\bold L|=k-1}\frac{(-1)^{||\bold L||-1}}{\bold L!}\bold
s^{\bold L}.$$ In particular, for any fixed values of $\bold s$,
$G(\bold s,\bold t)$ is a $\tau$-function for the KdV hierarchy.
\end{theorem}

At a final remark, it would be interesting to prove that
$\alpha_{\bold L}$ in Theorem 3.3 are positive for all $\bold L\in
N^\infty$. This problem is kindly pointed out to us by a referee.

More generally the question can be formulated as following: two
sequences $\alpha_{\bold L}$ and $\beta_{\bold L}$ with
$\alpha_{\bold 0}=\beta_{\bold 0}=1$ are said to be inverse to each
other if they satisfy
$$\left(\sum_{\bold L}\alpha_{\bold L}\bold s^{\bold L}\right)\cdot\left(\sum_{\bold L}\beta_{\bold L}\bold s^{\bold L}\right)=1.$$
Find sufficient conditions on $\beta_{\bold L}$ such that
$\alpha_{\bold L}>0$ for all $\bold L$.

We conjecture that $\alpha_{\bold L}$ are positive when $\sum_{\bold
L}\beta_{\bold L}\bold s^{\bold L}$ equals any of the following.
$$\sum_{\bold
L}\frac{(-1)^{||\bold L||}}{\bold L!(2|\bold{L}|+1)!!}\bold s^{\bold
L},\quad \sum_{\bold L}\frac{(-1)^{||\bold L||}}{\bold
L!(2|\bold{L}|-1)!!}\bold s^{\bold L},\quad \sum_{\bold
L}\frac{(-1)^{||\bold L||}}{\bold L!|\bold{L}|!}\bold s^{\bold L}.$$
The latter two arise when we consider Hodge integrals involving
$\lambda$ classes \cite{LX}.

For works on the positivity criteria of coefficients of reciprocal
power series of a single variable, see for example \cite{La}.
However it seems there is no literature dealing with the
coefficients of reciprocal series of several variables.

\end{document}